\documentclass[12pt]{article}

\usepackage{latexsym,amssymb,upref,amsmath,amsthm,amsfonts}
\usepackage {color}
\usepackage {pstricks, pst-node}
\usepackage {graphicx}
\usepackage {latexsym}
\usepackage [T1]{fontenc}
\usepackage [cp1250]{inputenc}
\usepackage{bbm}
\usepackage{bbold}

\newcommand{\achr}{\mathrm{achr}}
\newcommand{\ro}{\mathbb R}
\newcommand{\co}{\mathbb C}
\newcommand{\bp}{\noindent\textit{Proof. }}
\newcommand{\ep}{\blacksquare}
\newcommand{\exc}{\mathrm{exc}}
\newcommand{\cM}{\mathcal{M}}

\newcommand\cQ{\mathcal{Q}}

\newtheorem{theorem}{Theorem}

\newtheorem{proposition}[theorem]{Proposition}

\newtheorem{claim}{Claim}

\title{The achromatic number of $K_6\square K_7$ is $18$}

\author{Mirko Hor\v n\' ak\thanks{\noindent E-mail address: mirko.hornak$@$upjs.sk}\\
	Institute of Mathematics, P.J. \v{S}af\'arik University\\
	Jesenn\'a 5, 040\ 01 Ko\v{s}ice, Slovakia}

\begin{document}
\maketitle

\begin{abstract}
A vertex colouring $f:V(G)\to C$ of a graph $G$ is complete if for any two distinct colours $c_1,c_2\in C$ there is an edge $\{v_1,v_2\}\in E(G)$ such that $f(v_i)=c_i$, $i=1,2$. The achromatic number of $G$ is the maximum number $\achr(G)$ of colours in a proper complete vertex colouring of $G$. In the paper it is proved that $\achr(K_6\square K_7)=18$. This result finalises the determination of $\achr(K_6\square K_q)$.
\end{abstract}


\section{Introduction}

Consider a finite simple graph $G$ and a finite colour set $C$. A vertex colouring $f:V(G)\to C$ is \textit{complete} if for any two distinct colours $c_1,c_2\in C$ one can find an edge $\{v_1,v_2\}\in E(G)$ ($\{v_1,v_2\}$ is usually shortened to $v_1v_2$). The \textit{achromatic number} of $G$, in symbols $\achr(G)$, is the maximum cardinality of $C$ admitting a proper complete vertex colouring of $G$. The invariant was introduced by Harary, Hedetniemi, and Prins in \cite{HaHePr}, where the following interpolation theorem was proved:

\begin{theorem}\label{ip}
	If $G$ is a graph, and $\chi(G)\le k\le\achr(G)$ for an integer $k$, there exists a $k$-element colour set $C$ and a proper complete vertex colouring $f:V(G)\to C$. $\qed$
\end{theorem}

Let $G\square H$ denote the Cartesian product of graphs $G$ and $H$ (the notation follows the monograph \cite{IK} by Imrich and Klav\v zar). So, $V(K_p\square K_q)=V(K_p)\times V(K_q)$, and $E(K_p\square K_q)$ consists of all edges $(x,y_1)(x,y_2)$ with  $y_1\ne y_2$ and all edges $(x_1,y)(x_2,y)$ with $x_1\ne x_2$. The problem of determining $\achr(K_p\square K_q)$ is motivated by the fact that, according to Chiang and Fu \cite{ChiF1}, $\achr(G\square H)\ge\achr(K_p\square K_q)$ for arbitrary graphs $G,H$ with $\achr(G)=p$ and $\achr(H)=q$. Clearly, $\achr(K_q\square K_p)=\achr(K_p\square K_q)$, so it is natural to suppose $p\le q$. The case $p\in\{2,3,4\}$ was solved by Hor\v n\'ak and Puntig\'an in \cite{HoPu}, and that of $p=5$ by Hor\v n\'ak and P\v cola in \cite{HoPc1,HoPc2}. 

\begin{proposition}[see \cite{B}]\label{bou}
	$\achr(K_6\square K_6)=18$.
\end{proposition}

More generally, in~\cite{ChiF2} Chiang and Fu proved that if $r$ is an odd projective plane order, then $\achr(K_{(r^2+r)/2}\square K_{(r^2+r)/2})=(r^3+r^2)/2$.

The aim of the present paper is to finalise the determination of $K_6\square K_q$ (for $q\ge8$ see Hor\v n\'ak \cite{Ho1,Ho2}) by proving

\begin{theorem}\label{main}
$\achr(K_6\square K_7)=18$.
\end{theorem}

For $k,l\in\mathbb Z$ we denote \textit{integer intervals} by
\begin{equation*}
[k,l]=\{z\in\mathbb Z:k\le z\le l\},\qquad [k,\infty)=\{z\in\mathbb Z:k\le z\}.
\end{equation*}
Further, for a set $A$ and $m\in[0,\infty)$ let $\binom Am$ be the set of $m$-element subsets of $A$.

If we suppose that $V(K_r)=[1,r]$ for $r\in[1,\infty)$, then a vertex colouring $f:V(K_p\square K_q)\to C$ can be conveniently desribed using the $p\times q$ matrix $M=M(f)$, in which the entry in the $i$th row and the $j$th column is $(M)_{i,j}=f(i,j)$. 

If $f$ is proper, then each line (row or column) of $M$ consists of pairwise distinct entries. 

If $f$ is complete, then each pair $\{c_1,c_2\}\in\binom C2$ (of colours in $C$) is \textit{good} in $M$ in the sense that either there are $i\in[1,p]$ and $j_1,j_2\in[1,q]$ such that $\{c_1,c_2\}=\{(M)_{i,j_1},(M)_{i,j_2}\}$, which is expressed by saying that the pair $\{c_1,c_2\}$ is \textit{row-based} (in $M$) or there are $i_1,i_2\in[1,p]$ and $j\in[1,q]$ such that $\{c_1,c_2\}=\{(M)_{i_1,j},(M)_{i_2,j}\}$, \textit{i.e.}, the pair $\{c_1,c_2\}$ is \textit{column-based} (in $M$). For $(i,j)\in[1,p]\times[1,q]$ and $C'\subseteq C\setminus\{(M)_{i,j}\}$ let $g(i,j,C')$ denote the number of pairs $\{(M)_{i,j},\gamma\}$ with $\gamma\in C'$ that are good in $M$, while this fact is documented by the $(i,j)$-copy of the colour $(M)_{i,j}$.

Let $\cM(p,q,C)$ be the set of $p\times q$ matrices $M$ with entries from $C$ such that entries of any line of $M$ are pairwise distinct, and all pairs in $\binom C2$ are good in $M$. Thus if $f:[1,p]\times[1,q]\to C$ is a proper complete vertex colouring of $K_p\square K_q$, then $M(f)\in\cM(p,q,C)$. 

Conversely, if $M\in\cM(p,q,C)$, it is immediate to see that the mapping $f_M:[1,p]\times[1,q]\to C$ determined by $f_M(i,j)=(M)_{i,j}$ is a proper complete vertex colouring of $K_p\square K_q$.

So, we have just proved

\begin{proposition}\label{general}
	If $p,q\in[1,\infty)$ and $C$ is a finite set, then the following statements are equivalent:
	
	$\mathrm(1)$ There is a proper complete vertex colouring of $K_p\square K_q$ using as colours elements of $C$.
	
	$\mathrm(2)$ $\mathcal M(p,q,C)\ne\emptyset$.\quad$\qed$
\end{proposition}

The following straightforward proposition comes from~\cite{Ho1}.

\begin{proposition}\label{wlog}
	If $p,q\in[1,\infty)$, $C,D$ are finite sets, $M\in\mathcal M(p,q,C)$, mappings $\rho:[1,p]\to[1,p]$, $\sigma:[1,q]\to[1,q]$, $\pi:C\to D$ are bijections, and $M_{\rho,\sigma}$, $M_{\pi}$ are $p\times q$ matrices defined by $(M_{\rho,\sigma})_{i,j}=(M)_{\rho(i),\sigma(j)}$ and $(M_{\pi})_{i,j}=\pi((M)_{i,j})$, then $M_{\rho,\sigma}\in\mathcal M(p,q,C)$ and $M_{\pi}\in\mathcal M(p,q,D)$.\quad$\qed$
\end{proposition}

Let $M\in\cM(p,q,C)$. The \textit{frequency} of a colour $\gamma\in C$ is the number of entries of $M$ equal to $\gamma$. An \textit{l-colour} (\textit {$l+$colour}) is a colour of frequency $l$ (at least $l$), and $C_l$ ($C_{l+}$) is the set of $l$-colours ($l+$colours). Further, for $k\in\{l,l+\}$ let $c_k=|C_k|$,
\begin{alignat*}{4}
\ro(i)&=\{(M)_{i,j}:j\in[1,q]\},\ &&\ro_k(i)=C_k\cap\ro(i),\ &&r_k(i)=|\ro_k(i)|,\ &&i\in[1,p],\\
\co(j)&=\{(M)_{i,j}:i\in[1,p]\},\ &&\co_k(j)=C_k\cap\co(j),\ &&c_k(j)=|\co_k(j),\ &&j\in[1,q].
\end{alignat*}
Finally, denote 
\begin{alignat*}{3}
\ro_2(i_1,i_2)&=C_2\cap\ro(i_1)\cap\ro(i_2),\ &&r_2(i_1,i_2)=|\ro_2(i_1,i_2)|,\ &&i_1,i_2\in[1,p],i_1\ne i_2,\\
\co_2(j_1,j_2)&=C_2\cap\co(j_1)\cap\co(j_2),\ &&c_2(j_1,j_2)=|\co_2(j_1,j_2)|,\ &&j_1,j_2\in[1,q],j_1\ne j_2.
\end{alignat*}

Considering a nonempty set $S\subseteq[1,p]\times[1,q]$ we say that a colour $\gamma\in C$ \textit{occupies a position in $S$} (\textit{appears in $S$}, \textit{has a copy in $S$} or simply \textit{is in $S$}) if there is $(i,j)\in S$ such that $(M)_{i,j}=\gamma$.

\section{Properties of a counterexample to Theorem~\ref{main}}

We prove Theorem~\ref{main} by the way of contradiction. It is well known that $\achr(G)\ge\achr(H)$ if $H$ is an induced subgraph of a graph $G$. So, by Proposition~\ref{bou}, $\achr(K_6\square K_7)\ge\achr(K_6\square K_6)=18$. Provided that Theorem~\ref{main} is false, by Theorem~\ref{ip} and Proposition~\ref{general} there is a set $C$ with $|C|=19$ and a matrix $M\in\cM(6,7,C)$; henceforth the whole notation corresponds to this (hypothetical) matrix $M$. 

For a colour $\gamma\in C$ denote $V_{\gamma}=f_M^{-1}(\gamma)\subseteq[1,6]\times[1,7]$, and let $N(V_{\gamma})$ be the neighbourhood of $V_{\gamma}$ (the union of neighborhoods of vertices in $V_{\gamma}$). The \textit{excess} of $\gamma$ is defined to be the maximum number $\exc(\gamma)$ of vertices in a set $S\subseteq N(V_{\gamma})$ such that the partial vertex colouring of $K_6\square K_7$, obtained by removing colours of $S$, is complete concerning the colour class $\gamma$.

\begin{claim}\label{exc}
	If $\gamma\in C_l$, then $\exc(\gamma)=-l^2+12l-18$.
\end{claim}

\bp The vertex colouring $f_M$ of $K_6\square K_7$ is proper, hence $|N(V_{\gamma})|=7l+l(6-l)-l=l(12-l)$. Further, $f_M$ is complete, and so each colour of $C\setminus\{\gamma\}$ occupies a position in $N(V_{\gamma})$. As a consequence, $\exc(\gamma)=l(12-l)-(19-1)=-l^2+12l-18$. $\ep$

\begin{claim}\label{eleven}
	The following statements  are true:
	
	$1.$ $c_1=0$;
	
	$2.$ if $l\in[7,\infty)$, then $c_l=0$;
	
	$3.$ $c_2\in[15,18]$;
	
	$4.$ $c_{3+}\in[1,4]$;
	
	$5.$ $c_{4+}\le c_2-15$;
	
	$6.$ if $c_{4+}=0$, then $c_{3+}=c_3=4$;
	
	$7.$ if $c_{4+}\ge1$, then $c_{3+}\le3$;
	
	$8.$ $c_{3+}+c_{4+}\le4$;
	
	$9.$ if $c_{5+}\ge1$, then $c_{3+}+c_{4+}\le3$;
	
	$10.$ $\Sigma=\sum_{i=3}^6ic_i\in[6,12]$;
	
	$11.$ if $i_1,i_2\in[1,6]$, $i_1\ne i_2$, then $r_2(i_1,i_2)\le3$.
\end{claim}

\bp 1. If $\gamma\in C_1$, then, by Claim~\ref{exc}, $\exc(\gamma)=-7<0$, a contradiction.

2. If $\gamma\in C_l$ for some $l\in[7,\infty)$, by the pigeonhole principle the colouring $f_M$ is not proper, a contradiction.

3. By Claims~\ref{eleven}.1, \ref{eleven}.2 we have $c_2+c_{3+}=\sum_{i=2}^6c_i=|C|=19$ and $\sum_{i=2}^6ic_i=|V(K_6\square K_7)|=42$, hence $2c_2+6(19-c_2)=2c_2+6c_{3+}\ge\sum_{i=2}^6ic_i\ge2c_2+3c_{3+}=2c_2+3(19-c_2)$, which yields $114-4c_2\ge42\ge57-c_2$ and $15\le c_2\le18$.

4. A consequence of Claim~\ref{eleven}.3.

5. We have $3\cdot19-c_2+c_{4+}=3(c_2+c_3+c_{4+})-c_2+c_{4+}\le\sum_{i=2}^6 ic_i=42$ and $c_{4+}\le c_2-15$.

6. If $c_{4+}=0$, then $c_2+c_3=19$, $2c_2+3c_3=42$ and $c_{3+}=c_3=42-2\cdot19=4$.

7. The assumption $c_{4+}\ge1$ and $c_{3+}=4$ would mean $\Sigma\ge3\cdot3+4=13>12$, a contradiction.

8. If $c_{4+}=0$, then $c_{3+}+c_{4+}=c_{3+}\le4$. With $c_{4+}=1$ we have, by Claim~\ref{eleven}.5, $1=c_{4+}\le c_2-15$, $c_2\ge16$, $c_3=19-c_2-c_{4+}\le2$, $c_{3+}\le3$ and $c_{3+}+c_{4+}\le4$. Finally, from $c_{4+}\ge2$ it follows that $2\le c_{4+}\le c_2-15$, $c_2\ge17$, $19=c_2+c_3+c_{4+}\ge17+c_3+2=19+c_3\ge19$, $c_2=17$, $c_3=0$, $c_{4+}=2=c_{3+}$ and $c_{3+}+c_{4+}=4$.

9. We have $2(19-c_{5+})+5c_{5+}=2(c_2+c_3+c_4)+5(c_5+c_6)\le\sum_{i=2}^6ic_i=42$, $3c_{5+}\le4$ and $c_3+2c_4=42-2(19-c_5-c_6)-5c_5-6c_6\le4-3c_{5+}$, hence $c_{5+}\ge1$ yields $c_{5+}=1$, $c_3+2c_4\le1$ and $c_{3+}+c_{4+}\le1+2=3$.

10. The assertion of Claim~\ref{eleven}.3 leads to $\sum_{i=3}^6 ic_i=\sum_{i=2}^6 ic_i-2c_2=42-2c_2\in[6,12]$.

11. The inequality is trivial if $r_2(i_1,i_2)=0$. If $\gamma\in\ro_2(i_1,i_2)$, then each colour of $\ro_2(i_1,i_2)\setminus\{\gamma\}$ contributes one to the excess of $\gamma$, hence, by Claim~\ref{exc}, $r_2(i_1,i_2)-1\le\exc(\gamma)=2$ and $r_2(i_1,i_2)\le3$. $\ep$
\vskip1mm

A set $D\subseteq C_2$ is of the \textit{type} 
$(m_1^{a_1}\negthinspace\dots m_k^{a_k},n_1^{b_1}\negthinspace\dots n_l^{b_l})$ if $(m_1,\dots,m_k)$, $(n_1,\dots,n_l)$ are decreasing sequences of integers from the interval $[1,|D|]$, $|\{i\in[1,6]:|D\cap\ro(i)|=m_i\}|=a_i\ge1$ for each $i\in[1,k]$, and  $|\{j\in[1,7]:|D\cap\co(j)|=n_j\}|=b_j\ge1$ for each $j\in[1,l]$; clearly, $\sum_{i=1}^k m_ia_i=2|D|=\sum_{j=1}^l n_jb_j$.

Forthcoming Claims \ref{c23}, \ref{r23}, \ref{c22}, \ref{r22}, and \ref{3C2} state that certain types of 2- and 3-element subsets of $C_2$ are impossible in a matrix $M$ contradicting Theorem~\ref{main}. The mentioned claims are proved by contradiction. When proving that $M$ avoids a type $T$, we suppose that there is $D\subseteq C_2$, which is of the type $T$ (without explicitly mentioning it), and we reach a contradiction with some of the properties following from the fact that $M\in\cM(6,7,C)$.

\begin{claim}\label{c23}
No set $\{\alpha,\beta\}\subseteq C_2$ is of the type $(1^4,2^2)$.
\end{claim}

\bp Since we have at our disposal Proposition~\ref{wlog}, we may suppose without loss of generality that $(M)_{1,1}=\alpha=(M)_{3,2}$ and $(M)_{2,1}=\beta=(M)_{4,2}$. We use (w) to express briefly that it is Proposition~\ref{wlog}, which enables us to simplify our reasoning by restricting our attention to matrices with a special structure. With $A=\co(1)\cup\co(2)$ we have $|A|\le10$. If $\gamma\in C\setminus A$, the fact that both pairs $\{\gamma,\alpha\}$ and $\{\gamma,\beta\}$ are good forces $\gamma$ to occupy a position in $S_{\alpha}=\{1,3\}\times[3,7]$ and in $S_{\beta}=\{2,4\}\times[3,7]$ as well. So, $|C\setminus A|\le10$ and $|C\setminus A|=|C|-|A|\ge9$. By Claim~\ref{eleven}.4,
$|(C\setminus A)\cap C_2|=|C\setminus A|-|(C\setminus A)\cap C_{3+}|\ge9-c_{3+}\ge5$, hence there is $(i,k)\in\{1,3\}\times\{2,4\}$ such that $\delta\in\ro_2(i,k)\cap(C\setminus A)$. If $(i,k)=(1,2)$, (w) $(M)_{1,3}=\delta=(M)_{2,4}$.

If $|C\setminus A|=10$, then all ten positions in both $S_{\alpha}$, $S_{\beta}$ are occupied by colours of $C\setminus A$, and all twelve bullet positions in Fig.~1 are occupied by colours of $(C\setminus  A)\setminus\{\delta\}$, which means that $\exc(\delta)\ge12-|(C\setminus A)\setminus\{\delta\}|=12-(10-1)=3$ in contradiction to Claim~\ref{exc}.

If $|C\setminus A|=9$, then $\co(1)\cap\co(2)=\{\alpha,\beta\}$. Let $B$ be the set of four colours occupying a position in $[5,6]\times[1,2]$. Using $\exc(\alpha)=\exc(\beta)=2$ we see that exactly two positions in $[1,4]\times[3,7]$ are occupied by a colour of $B$. Thus, if $\varepsilon\in B$ is not in $[1,4]\times[3,7]$, it must be in $[5,6]\times[3,7]$, and so $\varepsilon\in C_2$. Then, however, the number of pairs $\{\zeta,\varepsilon\}$ with $\zeta\in(C\setminus A)\cap C_2$ that are good is at most four, while $|(C\setminus A)\cap C_2|\ge|C\setminus A|-|C_{3+}|\ge9-4=5$, a contradiction. 

Provided that $(i,k)\ne(1,2)$, a contradiction can be reached in a similar manner. $\ep$

\[
\begin{pmatrix}\label{col}
\alpha &. &\delta &\bullet &\bullet &\bullet &\bullet\\
\beta &. &\bullet &\delta &\bullet &\bullet &\bullet\\
. &\alpha &\bullet &\bullet &. &. &.\\
. &\beta &\bullet &\bullet &. &. &.\\
. &. &. &. &. &. &.\\
. &. &. &. &. &. &.\\
\end{pmatrix}  
\begin{pmatrix}
\alpha &\beta &. &. &. &. &.\\
. &. &\alpha &\beta &. &. &.\\
\gamma &\bullet &\bullet &\bullet &. &. &.\\
\bullet &\gamma &\bullet &\bullet &. &. &.\\
\bullet &\bullet &. &. &. &. &.\\
\bullet &\bullet &. &. &. &. &.
\end{pmatrix}
\begin{pmatrix}
\alpha &\beta &. &. &\circ &. &.\\
\zeta/. &./\zeta &\alpha &\beta &\circ &\circ &\circ\\
\gamma &. &. &. &. &\bullet &\bullet\\
. &\gamma &. &. &. &\bullet &\bullet\\
. &. &. &. &\zeta &. &.\\
. &. &. &. &. &. &.
\end{pmatrix}
\]
\[
\hskip0.5mm\mathrm{Fig.\ 1}\hskip35mm\mathrm{Fig.\ 2}\hskip35mm\mathrm{Fig.\ 3}\ 
\]

\begin{claim}\label{c22}
	No set $\{\alpha,\beta\}\subseteq C_2$ is of the type $(2^11^2,2^2)$.
\end{claim}

\bp Now (w) $(M)_{1,1}=\alpha=(M)_{2,2}$ and $(M)_{2,1}=\beta=(M)_{3,2}$. With $A=\co(1)\cup\co(2)\cup\ro(2)$ each colour $\gamma\in  C\setminus A$ has a copy in $\{i\}\times[3,7]$, $i=1,3$ ($\{\gamma,\alpha\}$ and $\{\gamma,\beta\}$ are good). From $|A|\le15$ it follows that $|C\setminus A|\ge19-15=4$, and then $C\setminus A\subseteq C_{3+}$: indeed, if $\delta\in(C\setminus A)\cap C_2$, then $\exc(\delta)\ge|(C\setminus A)\setminus\{\delta\}|\ge3$, a contradiction. Thus $C_2\subseteq A$, $c_2\le15$, hence, by Claims~\ref{eleven}.3, 2.4, $c_2=15$, $c_3=c_{3+}=4$, $C\setminus A=C_{3+}=C_3$, each colour of $C_2\setminus\{\alpha,\beta\}$ has exactly one copy in  $([1,6]\times[1,2])\cup(\{2\}\times[3,7])$, and (w) $\varepsilon=(M)_{1,2}$, $\zeta=(M)_{3,1}$ are (distinct) 2-colours.

First note that $\varepsilon,\zeta\notin\ro_2(1,3)$, for otherwise $\max(\exc(\varepsilon),\exc(\zeta))\ge c_3=4$. So, the second copies of $\varepsilon,\zeta$ are in $[4,6]\times[3,7]$, and the pair $\{\varepsilon,\zeta\}$ is good in the corresponding $3\times5$ submatrix of $M$. 

If the pair $\{\varepsilon,\zeta\}$ is column-based, (w) $\varepsilon=(M)_{4,3}$ and $\zeta=(M)_{5,3}$, then, with a (2-) colour $\eta$ appearing in $\{2\}\times[4,7]$, both pairs $\{\eta,\varepsilon\}$ and $\{\eta,\zeta\}$ are good only if $\eta$ occupies a position in $\{1,3,6\}\times\{3\}$, a contradiction.

If the pair $\{\varepsilon,\zeta\}$ is row-based, (w) $\varepsilon=(M)_{4,3}$ and $\zeta=(M)_{4,4}$, consider six positions in $[5,6]\times[5,7]$. Since $r_3(1)=r_3(3)=4=c_3$, at most four of those positions are occupied by 3-colours and at least two of them by 2-colours. Let $B$ be the set of 2-colours having a copy in $([5,6]\times[5,7])\cup(\{2\}\times[5,7])$. If $\vartheta\in B$, then, having in mind that both pairs $\{\vartheta,\varepsilon\}$ and $\{\vartheta,\zeta\}$
are good, $\vartheta$ must have a copy in $\{(1,4),(3,3),(4,1),(4,2)\}$; this yields a contradiction, because $|B|\ge5$. $\ep$
 
 \begin{claim}\label{ic2}
 	If $j,l\in[1,7]$, $j\neq l$, then $c(j,l)\le2$. 
 \end{claim}
  
 \bp The assumption $c(j,l)\ge3$ would contradict Claim~\ref{c23} or Claim~\ref{c22}. $\ep$

\begin{claim}\label{r23}
	No set $\{\alpha,\beta\}\subseteq C_2$ is of the type $(2^2,1^4)$.
\end{claim}

\bp Here (w) $(M)_{1,1}=\alpha=(M)_{2,3}$ and $(M)_{1,2}=\beta=(M)_{2,4}$. With $A=\ro(1)\cup\ro(2)$ we have $|A|\le12$, each colour of $C\setminus A$ is in both sets $S_{\alpha}=[3,6]\times\{1,3\}$, $S_{\beta}=[3,6]\times\{2,4\}$, and $7\le|C\setminus A|\le8$. As $|(C\setminus A)\cap C_2|\ge3$, there is $(j,l)\in\{1,3\}\times\{2,4\}$ such that $\gamma\in(C\setminus A)\cap\co_2(j,l)$. 

If $(j,l)=(1,2)$, (w) $(M)_{3,1}=\gamma=(M)_{4,2}$. 

If $|C\setminus A|=8$, all eight positions in both sets $S_{\alpha}$, $S_{\beta}$ are occupied by colours of $C\setminus A$. Further, all ten bullet positions in Fig.~2, which are positions of vertices in $(N(V_{\alpha})\cup N(V_{\beta}))\cap N(V_{\gamma})$, are occupied by colours of $(C\setminus A)\setminus\{\gamma\})$, and so $\exc(\gamma)\ge10-(8-1)=3$, a contradiction. 

Suppose that $|C\setminus A|=7$ (and $A|=12$). For $m\in\{2,3+,4+\}$ and $n\in[0,2]$ let $C_m^n$ be the set of colours in $C_m$ having $n$ copies in $[5,6]\times[5,7]$, and let $c_m^n=|C_m^n|$. If $\delta\in C_2^1\cup C_3^2$, then, since the pairs $\{\delta,\alpha\}$, $\{\delta,\beta\}$ and $\{\delta,\gamma\}$ are good, $\delta$ must appear in $\{2\}\times[1,2]$, and so $c_2^1+c_3^2\le2$; further, $c_2^2=0$. Using Claim~\ref{eleven}.8 we obtain 
\begin{align}\label{=6}
6&=c_2^1+c_3^1+c_{4+}^1+2c_3^2+2c_{4+}^2\le c_2^1+c_3^2+\sum_{n=0}^2(c_3^n+c_{4+}^n)+c_{4+}^2\nonumber\\
&\le c_2^1+c_3^2+c_{3+}+c_{4+}\le2+4=6,
\end{align}
which implies
\begin{align}
c_3^0=c_{4+}^0=c_{4+}^1&=0,\label{=0}\\
c_{4+}&=c_{4+}^2\label{4+},\\
c_2^1+c_3^2&=2\label{=2},
\end{align}
$c_{3+}+c_{4+}=4$, and so, by Claim~\ref{eleven}.9, $c_{5+}=0$.

For $\delta\in\{\alpha,\beta\}$ choose a set $S_{\delta'}\subseteq S_{\delta}$ with $|S_{\delta'}|=7$ occupied by seven distinct colours of $C\setminus A$, and let $P=([3,6]\times[1,4])\setminus(S_{\alpha'}\cup S_{\beta'})$; then $|P|=2$. Since $|N(V_{\gamma})\cap([3,6]\times[1,4])|=10$, we have $2=\exc(\gamma)\ge10-(|P|+|(C\setminus A)\setminus\{\gamma\}|=4-|P|=2$, hence both positions in $P$ are necessarily occupied by a colour of $A$, and all sets $S_{\alpha'},S_{\beta'},P$ are unique. We express this property of $\gamma$ by saying that $\gamma$ is \textit{$A$-exact}. Besides that, the two positions in $P$ are occupied by two distinct colours of $A$, say $\lambda$ and $\mu$; indeed, otherwise the colour of $A$, which occupies both positions in $P$, by (\ref{4+}), would be a 5+colour, a contradiction. Let $P=\{(i_{\lambda},j_{\lambda}),(i_{\mu},j_{\mu})\}$, where $\lambda=(M)_{i_{\lambda},j_{\lambda}}$ and $\mu=(M)_{i_{\mu},j_{\mu}}$. The excess of both $\alpha,\beta$ is 2, therefore $(j_{\lambda},j_{\mu})\in\{1,3\}\times\{2,4\}$ (a colour occupying a position in $P$ contributes to the excess of either $\alpha$ or $\beta$, and $\alpha,\beta$ are contributing to the excess of each other).

The above reasoning concerning $\gamma$ can be repeated to prove that any colour in $(C\setminus A)\cap C_2$ is $A$-exact.

Suppose that $\varepsilon\in(C\setminus A)\cap C_2$; then $\varepsilon$ is $A$-exact and $\varepsilon\in\co_2(j',l')$, where $(j',l')\in\{1,3\}\times\{2,4\}$. Let $\{l_{\lambda},l_{\mu}\}=[1,4]\setminus\{j_{\lambda},j_{\mu}\}$. 

Assume first that $i_{\lambda}=i_{\mu}$. By Claim~\ref{c22}, $|(C\setminus A)\cap\co_2(j_{\lambda},j_{\mu})|\le2$. If $(j',l')\ne(j_{\lambda},j_{\mu})$, then either $\varepsilon=(M)_{i_{\lambda},l_{\lambda}}$ or $\varepsilon=(M)_{i_{\lambda},l_{\mu}}$.

The second possibility is $i_{\lambda}\ne i_{\mu}$. By Claim~\ref{c22}, $|(C\setminus A)\cap\co_2(l_{\lambda},l_{\mu})|\le2$. On the other hand, if $(j',l')\ne(l_{\lambda},l_{\mu})$, then either $\varepsilon=(M)_{i_{\lambda},j_{\mu}}$ or $\varepsilon=(M)_{i_{\mu},j_{\lambda}}$.

In both cases $|(C\setminus A)\cap C_2|\le2+2=4$ and 
\begin{equation}\label{ge3}
|(C\setminus A)\cap C_{3+}|\ge7-4=3. 
\end{equation}
From (\ref{=6}) and (\ref{4+}) we obtain $(C\setminus A)\cap C_{3+}\subseteq C_3^1\cup C_{4+}^2$, hence $c_3^1+c_{4+}^2\ge3$. Let us show the following:
\vskip1mm

\noindent (*) No 3+colour occupies a position in $[3,4]\times[5,7]$, and $c_2^1\ge1$.
\vskip1mm

Because of (\ref{=0}) and (\ref{4+}) we know that colours of $(C\setminus A)\cap C_{3+}$ appear only in $([3,6]\times[1,4])\cup([5,6]\times[5,7])$. If $c_{4+}\ge1$, then, by Claim~\ref{eleven}.7, $3\ge c_{3+}\ge c_3^1+c_{4+}^2\ge3$, $c_{3+}=c_3^1+c_{4+}^2=3$, $C_{3+}=C_3^1\cup C_{4+}^2=(C\setminus A)\cap C_{3+}$, $c_3^2=0$, $c_2^1=2$ (see (\ref{=2})), and so (*) is true. 

If $c_{4+}=0$, from (\ref{=6}), (\ref{=2}) and Claim~\ref{eleven}.4 it follows that $c_3^1+c_3^2=4=c_{3+}$, $C_{3+}=C_3^1\cup C_3^2=((C\setminus A)\cap C_{3+})\cup C_3^2$ and, by (\ref{ge3}), $c_3^1\ge3$; since a colour of $C_3^2$ is only in $(\{2\}\times[1,2])\cup([5,6]\times[5,7])$, $c_3^2\le1$ and $c_2^1\ge1$, (*) is true again.

Now, by (*), six positions in $[3,4]\times[5,7]$ are occupied by six distinct 2-colours belonging to $A\setminus\{\lambda,\mu\}$, and there is a colour $\zeta\in C_2^1$, (w) $\zeta=(M)_{5,5}$, see Fig.~3. If a 2-colour $\eta$ appears in a bullet position, then, since the pair $\{\eta,\zeta\}$ is good, the second copy of $\eta$ must occupy a circle position. In such a case, however, it is easy to check that there is a set $\{\vartheta,\iota\}\subseteq A\cap C_2$ of 2-colours occupying two bullet positions and two circle positions, which contradicts either Claim~\ref{c23} or Claim~\ref{c22}.

The case $(j,l)\ne(1,2)$ can be treated similarly. $\ep$

\begin{claim}\label{r22}
No set $\{\alpha,\beta\}\subseteq C_2$ is of the type $(2^2,2^11^2)$.
\end{claim}

\bp Let (w) $(M)_{1,1}=\alpha=(M)_{2,2}$ and $(M)_{1,2}=\beta=(M)_{2,3}$. First of all, we have $\ro_2(1,2)=\{\alpha,\beta\}$. Indeed, if (w) $\gamma\in\ro(1,2)\setminus\{\alpha,\beta\}$, then, by Claim~\ref{r23}, necessarily $(M)_{1,3}=\gamma=(M)_{2,1}$. Each colour $\delta\in C\setminus\ro_2(1,2)$ occupies at least two positions in $[3,6]\times[1,3]$ (all pairs $\{\delta,\alpha\}$, $\{\delta,\beta\}$, $\{\delta,\gamma\}$ are good), hence $|C|\le11+\lfloor\frac{4\cdot3}2\rfloor=17<19$, a contradiction.

With $A=\ro(1)\cup\ro(2)\cup\co(2)$ any colour $\gamma\in C\setminus A$ occupies a position in $[3,6]\times\{1\}$ as well as in $[3,6]\times\{3\}$, hence $|C\setminus A|\le4$, $|A|\le16$, $|C\setminus A|=19-|A|\ge3$ and $|A|\ge15$. 

Assume first that $|A|=15$ and $|C\setminus A|=4$, which yields $C\setminus A\subseteq C_{3+}$ (a 2-colour $\gamma\in C\setminus A$ would satisfy $\exc(\gamma)\ge3$), $c_{3+}=c_3=4$ and $A=C_2$. For colours $\gamma=(M)_{1,3}$ and $\delta=(M)_{2,1}$ their second copies appear in $[3,6]\times(\{2\}\cup[4,7])$, and 
the pair $\{\gamma,\delta\}$ is good in the corresponding $4\times5$ submatrix of $M$. However, at most one of $\gamma,\delta$ is in $[3,6]\times\{2\}$, hence $\{\gamma,\delta\}$ is good in the submatrix of $M$ corresponding to $[3,6]\times[4,7]$.

If the pair $\{\gamma,\delta\}$ is column-based, (w) $\gamma=(M)_{3,4}$ and $\delta=(M)_{4,4}$, at most one of colours in $[5,6]\times\{2\}$ belongs to $\ro(1)\cup\ro(2)$, hence there is a 2-colour $\varepsilon$ and $i\in[5,6]$ such that $\varepsilon=(M)_{i,2}=(M)_{11-i,4}$ (so that both pairs $\{\varepsilon,\gamma\}$, $\{\varepsilon,\delta\}$ are good). For every colour $\zeta\notin(\co(2)\cup\{(M)_{i,4}\}$ occupying a position in $[1,2]\times[5,7]$ (the number of such colours is at least 4) there is $\eta\in\{\gamma,\delta,\varepsilon\}$ such that the pair $\{\eta,\zeta\}$ is not good, a contradiction.

If the pair $\{\gamma,\delta\}$ is row-based, (w) $\gamma=(M)_{3,4}$ and $\delta=(M)_{3,5}$. If a colour $\varepsilon$ occupies a position in $[4,6]\times\{2\}$ and does not belong to $\ro(1)\cup\ro(2)$ (there are at least two such colours), it must appear in $\{3\}\times[6,7]$ (pairs $\{\varepsilon,\gamma\}$ and $\{\varepsilon,\delta\}$ are good), (w) $(M)_{4,2}=\varepsilon=(M)_{3,6}$ and $(M)_{5,2}=\zeta=(M)_{3,7}$. If a 2-colour $\eta$ is in $\{6\}\times[4,7]$, then $\eta=(M)_{3,2}$ (all pairs $\{\eta,\vartheta\}$ with $\vartheta\in\{\gamma,\delta,\varepsilon,\zeta\}$ are good), $r_3(6)\ge2+3=5>c_3$, and so the colouring $f_M$ is not proper, a contradiction.

From now on $|A|=16$ and $|C\setminus A|=3$. Suppose first that $C\setminus A\subseteq C_{3+}$. From $c_{3+}\le4$ we obtain $|A\cap C_{3+}|\le1$. 

If $(\ro(1)\cup\ro(2))\cap C_{3+}=\emptyset$, (w) $\gamma=(M)_{3,2}$, $\delta=(M)_{4,2}$, $\varepsilon=(M)_{5,2}$ are 2-colours, and their second copies appear in $[3,6]\times[4,7]$. Let $J=\{j\in[4,7]:\co(j)\cap\{\gamma,\delta,\varepsilon\}\ne\emptyset\}$; by Claim~\ref{ic2} we know that $2\le|J|\le3$. If $(i,j)\in S=\{(1,3),(2,1)\}\cup([1,2]\times([4,7]\setminus J))$, then $g(i,j,\{\gamma,\delta,\varepsilon\})=0$; note that $|S|=10-2|J|$. On the other hand, it is easy to see that the number of pairs $(i,j)\in[3,6]\times[4,7]$, satisfying $g(i,j,\{\gamma,\delta,\varepsilon\})=3$, is less than $|S|$ (at most 3 if $|J|=3$ and at most 4 if $|J|=2$). Thus, there is a 2-colour $\zeta$ in $S$ and $\eta\in\{\gamma,\delta,\varepsilon\}$ such that the pair $\{\zeta,\eta\}$ is not good. 

If $|(\ro(1)\cup\ro(2))\cap C_{3+}|=1$, then $c_3=c_{3+}=4$ and $c_2(2)=6$. 

Suppose first that both $\gamma=(M)_{1,3}$ and $\delta=(M)_{2,1}$ are 2-colours. The second copies of $\gamma,\delta$ are then in $[3,6]\times[4,7]$, for if not, $\max(\exc(\gamma),\exc(\delta))\ge1+|C\setminus A|=4$. 

If the pair $\{\gamma,\delta\}$ is column-based, (w) $\gamma=(M)_{3,4}$ and $\delta=(M)_{4,4}$, then $(M)_{5,2}=\varepsilon=(M)_{6,4}$ and $(M)_{6,2}=\zeta=(M)_{5,4}$ (all pairs $\{\varepsilon,\gamma\},\{\varepsilon,\delta\},\{\zeta,\gamma\},$ $\{\zeta,\delta\}$ are good). For $(i,j)\in[1,2]\times[5,7]$ then $g(i,j,\{\gamma,\delta,\varepsilon,\zeta\})=1$, and at least three positions in $[1,2]\times[5,7]$ are occupied by a 2-colour that is in $[3,6]\times[5,7]$; on the other hand, for $(i,j)\in[3,6]\times[5,7]$ we have $g(i,j,\{\gamma,\delta,\varepsilon,\zeta\})\le2$, a contradiction.

If the pair $\{\gamma,\delta\}$ is row-based, (w) $\gamma=(M)_{3,4}$ and $\delta=(M)_{3,5}$. Then $g(i,j,\{\gamma,\delta\})=0$ for $(i,j)\in[4,6]\times\{2\}$ and $g(i,j,\{\gamma,\delta\})\le1$ for $(i,j)\in[4,6]\times[4,7]$; this leads to a contradiction, since at least one of colours in $[4,6]\times\{2\}$ has its second copy in $[4,6]\times[4,7]$.

So, one of $\gamma,\delta$ is a 2-colour and the other a 3-colour, (w) $\gamma\in C_2$ and $\delta\in C_3$. As above, the second copy of $\gamma$ appears in $[3,6]\times[4,7]$, (w) $\gamma=(M)_{3,4}$. All colours of the set $B=\{\varepsilon,\zeta,\eta,\vartheta\}$, where $\varepsilon=(M)_{3,2}$, $\zeta=(M)_{4,2}$, $\eta=(M)_{5,2}$ and $\vartheta=(M)_{6,2}$, are 2-colours. By Claim~\ref{c23}, the second copy of a colour $\iota\in B$ does not appear in $\co(1)\cup\co(3)$, hence is in $[3,6]\times[4,7]$ and, additionally, in $\ro(3)\cup\co(4)$, provided that $\iota\ne\varepsilon$ (the pair $\{\iota,\gamma\}$ is good). Clearly, $|B\cap\ro(3)|\le3$, since otherwise $\exc(\varepsilon)\ge3$. So, by Claim~\ref{ic2}, with $B'=\{\zeta,\eta,\vartheta\}$ we have $1\le|B'\cap\co(4)|\le2$.

If $|B'\cap\co(4)|=2$, (w) $\eta=(M)_{4,4}$, $\zeta=(M)_{4,5}$ (here we use Claim~\ref{c22}) and $\vartheta=(M)_{3,5}$. For both $l\in[6,7]$ then $g(2,l,B'\cup\{\gamma\})=0$. This leads to a contradiction, since $(M)_{2,6},(M)_{2,7}$ are 2-colours, and $g(i,j,B'\cup\{\gamma\})=4$ only if $(i,j)=(6,4)$.

If $|B'\cap\co(4)|=1$, (w) $\zeta=(M)_{5,4},\eta=(M)_{3,5}$ and $\vartheta=(M)_{3,6}$, then $g(2,7,B'\cup\{\gamma\})=0$. A contradiction follows from the fact that $g(i,j,B'\cup\{\gamma\})\le3$ for each $(i,j)$.

Now suppose that $(C\setminus A)\cap C_2\ne\emptyset$, (w) $\gamma=(M)_{3,1}=(M)_{4,3}\in(C\setminus A)\cap C_2$. For $m\in\{2,3,3+\}$, $n\in[1,2]$ let $C_m^n$ be the set of colours in $C_m$  having $n$ copies in $[5,6]\times[4,7]$ and $c_m^n=|C_m^n|$; then

\begin{equation}\label{8}
c_2^1+c_{3+}^1+2c_{3+}^2=8.
\end{equation}
Since $g(i,j,\{\alpha,\beta,\gamma\})=0$ for $(i,j)\in[5,6]\times[4,7]$ and $g(i,j,\{\alpha,\beta,\gamma\})=3$ if and only if $(i,j)\in S=\{(1,3),(2,1),(3,2),(4,2)\}$, we have

\begin{equation}\label{le4}
c_2^1+c_3^2\le4.
\end{equation}

Let us first show that $c_2^1\le3$. Indeed, if $c_2^1=4$, it is easy to see that all pairs $\{\delta,\varepsilon\}\in\binom{C_2^1}2$ are good only if there is $i\in[5,6]$ such that $C_2^1\subseteq\ro(i)$. This immediately implies $c_{3+}^2=0$ and, by Claim~\ref{eleven}.4 and (\ref{8}), $4\ge c_{3+}\ge c_{3+}^1=4$, $c_{3+}=4$ and $C_{3+}\subseteq\ro(11-i)$. Then $\delta=(M)_{11-i,2}\in C_2$, the second copy of $\delta$ is in $[3,4]\times[4,7]$ (by Claim~\ref{c23}), hence at least one of pairs $\{\delta,\varepsilon\}$ with $\varepsilon\in C_2^1$ is not good, a contradiction.

Further, we prove that 
\begin{equation}\label{sum}
c_{3+}^1+c_{3+}^2=c_{3+},
\end{equation}
which is equivalent to
\begin{equation}\label{eq}
C_{3+}^1\cup C_{3+}^2=C_{3+}.
\end{equation}

If $c_{4+}\ge1$, Claim~\ref{eleven}.7 yields $c_{3+}\le3$. Because of (\ref{8}) we obtain $2(c_{3+}^1+c_{3+}^2)=8+c_{3+}^1-c_2^1$, $c_{3+}^1+c_{3+}^2=\frac12(8+c_{3+}^1-c_2^1)\ge\frac12(8-3)=\frac52$, $3\ge c_{3+}\ge c_{3+}^1+c_{3+}^2\ge3$, and so (\ref{sum}) is true under the assumption $c_{4+}\ge1$ (implying $c_{3+}=3$). 
 
On the other hand, $c_{4+}=0$ implies $c_{3+}=c_3=4$ (Claim~\ref{eleven}.6). In this case, using (\ref{8}) and (\ref{le4}), we see that $8=(c_2^1+c_{3+}^2)+(c_{3+}^1+c_{3+}^2)\le4+c_{3+}=8$, hence $c_{3+}=4=c_2^1+c_{3+}^2=c_{3+}^1+c_{3+}^2$,
and (\ref{sum}) holds again. 

Note that now necessarily $|(C\setminus A)\cap C_2|=1$. Indeed, $|(C\setminus A)\cap C_2|=3$ is impossible by Claim~\ref{ic2}, since in such a case $c_2(1,3)\ge|(C\setminus A)\cap C_2|=3$. Moreover, by Claims~\ref{c23} and \ref{c22}, the assumption $|(C\setminus A)\cap C_2|=2$ would mean that for the unique colour $\delta\in(C\setminus A)\cap C_{3+}$ there is $i\in[5,6]$ such that $(M)_{i,1}=\delta=(M)_{11-i,3}$; however, according to (\ref{eq}), $\delta$ has an exemplar in $[5,6]\times[4,7]$, and so the colouring $f_M$ is not proper, a contradiction.

Thus 
\begin{equation}\label{extra+}
|(C\setminus A)\cap C_{3+}|=2.
\end{equation}
Because of a reasoning analogous to that above we see that each colour of $(C\setminus A)\cap C_{3+}$ occupies exactly one position in $[5,6]\times\{1,3\}$,
\begin{equation}\label{extra++}
(C\setminus A)\cap C_{3+}=C_{3+}^1,
\end{equation}
and then, using (\ref{eq}),
\begin{equation}\label{ro56}
C_{3+}\subseteq\ro(5)\cap\ro(6).
\end{equation}

Now if $c_{4+}\ge1$ (and, consequently, $c_{3+}=3$, which we have seen already), then, by (\ref{sum}), (\ref{extra+}) and (\ref{extra++}), $c_{3+}^1=2$ and $c_{3+}^2=1$; since $c_2^1\le3$, in such a case $c_2^1+c_{3+}^1+2c_{3+}^2\le7$ in contradiction with (\ref{8}).

Therefore, in the rest of the proof of Claim~\ref{r22} we work with $c_{4+}=0$, $c_{3+}=4$, $c_{3+}^1=2$, $c_3^2=c_{3+}^2=2$ and $c_2^1=2$, see (\ref{8}), (\ref{sum}), (\ref{extra+}), (\ref{extra++}). Moreover, all positions in $S$ are occupied by colours of $C_2^1\cup C_3^2$. If $\delta=(M)_{i,j}\in C_2^1$ with $(i,j)\in\{(1,3),(2,1)\}$, then, because of (\ref{extra+}), (\ref{extra++}) and (\ref{ro56}), $\exc(\delta)\ge1+|(C\setminus A)\cap C_{3+}|=1+c_{3+}^1=3$, a contradiction.

Thus $\{(M)_{1,3},(M)_{2,1}\}\subseteq C_{3+}^2$, and for a 2-colour $\varepsilon$ occupying a position in $[5,6]\times\{1,3\}$ (there are two such colours), by (\ref{ro56}) we have $\exc(\varepsilon)\ge c_{3+}-1=3$, a contradiction again. $\ep$

\begin{claim}\label{ir2}
If $i,k\in[1,6]$, $i\neq k$, then $r(i,k)\le2$. 
\end{claim}

\bp The assumption $r(i,k)\ge3$ would be in contradiction with Claim~\ref{r23} or Claim~\ref{r22}. $\ep$

\begin{claim}\label{3C2}
No set $\{\alpha,\beta,\gamma\}\subseteq C_2$ is of the type $(3^12^11^1,3^12^11^1)$.
\end{claim}
\bp Having in mind Claim~\ref{c22} (or else Claim~\ref{r22}), assume (w) $(M)_{1,1}=\alpha=(M)_{2,2}$, $(M)_{1,2}=\beta=(M)_{2,1}$ and  $(M)_{1,3}=\gamma=(M)_{3,1}$. Let $A^1$ be the set of colours occupying a position in $(\{1\}\times[4,7])\cup([4,6]\times\{1\})\cup\{(2,3),(3,2)\}$, $A^n$ the set of colours in $(\{n\}\times[4,7])\cup([4,6]\times\{n\})$ for $n=2,3$, $A_m^n=C_m\cap A^n$ and $a_m^n=|A_m^n|$ for $m\in\{2,3+\}$, $n\in[1,3]$. Then 
\begin{eqnarray}
|A^1|=a_2^1+a_{3+}^1=9,\label{=9}\\
|A^2|=a_2^2+a_{3+}^2=7,\label{=7}\\
A_1\cap A_2=\emptyset,\label{0}
\end{eqnarray}
since otherwise $\exc(\alpha)\ge3$. Moreover, $|A^3|\le7$ and $A^2\subseteq A^3$ (each pair $\{\gamma,\eta\}$ with $\eta\in A_2$ is good), hence, by (\ref{=7}),
\begin{equation}\label{=}
A_2=A_3.
\end{equation}

Let us show that distinct colours $\delta=(M)_{2,3}$, $\varepsilon=(M)_{3,2}$ (a consequence of (\ref{=9})) satisfy
\begin{equation}\label{de}
\{\delta,\varepsilon\}\subseteq C_{3+}.
\end{equation}
Indeed, if $\eta=(M)_{i,5-i}\in\{\delta,\varepsilon\}\cap C_2$ for some $i\in[2,3]$ and (w) $\eta=(M)_{4,4}$, then all colours appearing in $(\{i\}\times[4,7])\cup([4,6]\times\{5-i\})\cup\{(5-i,4),(4,i)\}$ belong to $A^2\setminus\{\eta\}$, hence, by (\ref{=7}) and (\ref{=}), $\exc(\eta)\ge9-(7-1)=3$, a contradiction.

Further, with $\zeta=(M)_{3,3}$ we have
\begin{equation}\label{zeta}
\zeta\in C_{3+}\cap A^1.
\end{equation}
To see it realise first that, since the pair $\{\zeta,\alpha\}$ is good and $f_M$ is proper, we get $\zeta\notin A^2\cup\{\delta,\varepsilon\}$ and $\zeta\in A^1$. Moreover, $\zeta\in C_{3+}$, because the assumption $\zeta\in\ro_2(1)$ ($\zeta\in\co_2(1)$) contradicts Claim~\ref{r22} (Claim~\ref{c22}, respectively). 

By (\ref{=9}), (\ref{=7}) and Claim~\ref{eleven}.4, we have
\begin{equation}\label{geq}
a_2^1+a_2^2\ge(9+7)-c_{3+}\ge12.
\end{equation}
Further, by (\ref{=9}), (\ref{=7}) and (\ref{de})--(\ref{geq}),
\begin{eqnarray}
5\le a_2^1\le6,\label{a21}\\
6\le a_2^2\le7.\label{a22}
\end{eqnarray}

Consider a colour $\eta\in A_2^1\cap\ro(1)$ (from(\ref{a21}) we see that there are at least two such colours), (w) $\eta=(M)_{1,j}=(M)_{4,l}$. Then from $g(1,j,A_2^2)\le g(1,j,A^2)=2$ and $g(4,l,A_2^2)\le g(4,l,A^2)\le4$ it follows that $a_2^2\le2+4=6$, hence, by (\ref{=7}), (\ref{a21}), (\ref{a22}), $a_2^1=a_2^2=6$ and $a_{3+}^2=1$.

Suppose first $\zeta\in\co(1)$ so that all positions in $\{1\}\times[4,7]$ are occupied by colours of $A_2^1$. If $\eta=(M)_{1,j}$, $j\in[4,7]$, then proceeding similarly as above it is easy to see that both positions in $[2,3]\times\{j\}$ are occupied by colours of $A_2^2$. Thus $\ro_2(2,3)$ consists of at least two colours of $A_2^2$. By Claim~\ref{ir2} we obtain $r_2(2,3)=2$, (w) $(M)_{2,4}=\vartheta=(M)_{3,5}$ and $(M)_{2,5}=\iota=(M)_{3,4}$. Now $\kappa=(M)_{1,6}$ satisfies $\kappa\in\co(4)\cup\co(5)$ (the pair $\{\kappa,\vartheta\}$ is good) and, analogously, $\lambda=(M)_{1,7}\in\co(4)\cup\co(5)$. By Claim~\ref{r23}, the copies of $\kappa,\lambda$ that are in $[4,6]\times[4,5]$ do not share a row, (w) one of them is in $\ro(4)$ and the other in $\ro(5)$. Then, however, reasoning similarly as above again, all positions in $[4,5]\times[2,3]$ are occupied by colours of $A_2^2$. Consequently, the unique colour of $A_{3+}^2$ is $(M)_{6,2}=(M)_{6,3}$, and $f_M$ is not proper.

If $\zeta\in\ro(1)$, (w) $\zeta=(M)_{1,7}$, then all positions in $(\{1\}\times[4,6])\cup([4,6]\times\{1\})$ are occupied by colours of $A_2^1$, which implies that all positions in $([2,3]\times[4,6])\cup([4,6]\times[2,3])$ are occupied by colours of $A_2^2$. So, the unique colour of $A_{3+}^2$ is $(M)_{2,7}=(M)_{3,7}$, a contradiction. $\ep$
 
\section{Final analysis} 

We are now ready to do the final analysis. Suppose (w) that $r_2(1)\ge r_2(i)$ for $i\in[2,6]$, which, by Claim~\ref{eleven}.3, implies 
\begin{equation}\label{r2}
7\ge r_2(1)\ge\left\lceil\frac{2c_2}6\right\rceil\ge\left\lceil\frac{30}6\right\rceil=5.
\end{equation}
Given $r_2(1)$ we assume (w) that the sequence $S=(r_2(1,i))_{i=2}^6$ is nonincreasing. Since $r_2(1)\in[5,7]$, we have $r_2(1,2)\ge\lceil\frac{r_2(1)}5\rceil\ge1$, $r_2(1,6)\le\lfloor\frac{r_2(1)}{5}\rfloor=1$, and Corollary~\ref{ir2} yields $r_2(1,2)\le2$. We suppose (w) that 
\begin{equation*}
j\in[1,r_2(1)]\Rightarrow(M)_{1,j}\in C_2,
\end{equation*}
and, more precisely,
\begin{eqnarray*}
(M)_{1,1}=\alpha,\ (M)_{1,2}=\beta,\ (M)_{1,3}=\gamma,\ (M)_{1,4}=\delta,\ (M)_{1,5}=\varepsilon,\\
r_2(1)\ge6\Rightarrow(M)_{1,6}=\zeta,\qquad r_2(1)=7\Rightarrow(M)_{1,7}=\eta.
\end{eqnarray*}
Let $p$ be the smallest integer in $[2,6]$ such that $r_2(1,i)\le1$ for every $i\in[p,6]$; $p$ is correctly defined since $r_2(1,6)\le1$. Then 
\begin{equation*}
r_2(1,i)=1\Leftrightarrow i\in[p,r_2(1)+3-p],
\end{equation*}
and, counting the number of positions in $[2,6]\times[1,7]$ occupied by colours of $\ro_2(1)$, we obtain $2(p-2)\le r_2(1)\le2(p-2)+(7-p)$, which yields
\begin{equation}\label{p}
r_2(1)-3\le p\le\left\lfloor\frac{r_2(1)+4}2\right\rfloor\le5.
\end{equation}
Moreover, because of Claims~\ref{r23} and \ref{r22} we have
\begin{align*}
p\ge3&\Rightarrow((M)_{2,1}=\beta\land(M)_{2,2}=\alpha),\\
p\ge4&\Rightarrow((M)_{3,3}=\delta\land(M)_{3,4}=\gamma),\\
p=5&\Rightarrow((M)_{4,5}=\zeta\land(M)_{4,6}=\varepsilon).
\end{align*}

Let $q_j=|\ro_2(1)\cap\co(j)|$ for $j\in[1,7]$. By Claim~\ref{3C2} we know that a 2-colour $\mu$, which occupies a position in $\{1\}\times[2p-3,r_2(1)]$, satisfies $\mu\notin\co(j)$ for every $j\in[1,2p-4]$, hence $q_j=2$ for any $j\in[1,2p-4]$, and
\begin{equation}\label{sum}
\sum_{j=2p-3}^7q_j=2[r_2(1)-(2p-4)]=2r_2(1)+8-4p;
\end{equation}
clearly, 
\begin{equation}\label{qj}
j\in[2p-3,7]\Rightarrow q_j\le3,
\end{equation}
since with $q_j\ge4$ and $\mu\in\ro_2(1)\cap\co(j)$ for some $j\in[2p-3,7]$ we have $\exc(\mu)\ge q_j-1\ge3$.  Notice also that 
\begin{equation}\label{qj+}
j\in[2p-3,r_2(1)]\Rightarrow1\le q_j\le\min(3,r_2(1)+4-2p),
\end{equation}
because 2-colours occupying a position in $[1,6]\times\{j\}$ are distinct from $2p-4$ (2-)colours appearing in $\{1\}\times[1,2p-4]$. Moreover, we assume (w) that
the sequence $(q_j)_{j=2p-3}^{r_2(1)}$ is nonincreasing, and that, if $(r_2(1),p)=(5,2)$, the sequence $(q_1,q_2,q_3,q_4,q_5)$ is nonincreasing.

For every pair $(r_2(1),p)$ obeying (\ref{r2}) and (\ref{p}), we analyse the set $\cQ(r_2(1),p)$ of sequences $(q_j)_{j=2p-3}^7$ satisfying all restrictions (\ref{sum})--(\ref{qj+}). More precisely,  we show that the assumption that $M$ is characterised by an arbitrary sequence $Q\in \cQ(r_2(1),p)$ leads to a contradiction, mostly because of $\Sigma\ge13$ (a contradiction to Claim~\ref{eleven}.10) or the existence of a line of $M$ containing at least five 3+colours (by Claim~\ref{eleven}.4 then the colouring $f_M$ is not proper). 

The structure of the sets $\cQ(r_2(1),p)$ with $(r_2(1),p)\ne(5,2)$ follows: 
\begin{align*}
\cQ(7,5)=&\emptyset,\\
\cQ(7,4)=&\{(3,2,1),(2,2,2)\},\\ 
\cQ(6,5)=&\{(0)\},\\
\cQ(6,4)=&\{(2,2,0),(2,1,1),(1,1,2)\},\\
\cQ(6,3)=&\{(3,3,1,1,0),(3,2,2,1,0),(3,2,1,1,1),(3,1,1,1,2),(2,2,2,2,0),\\
&(2,2,2,1,1),(2,2,1,1,2),(2,1,1,1,3)\},\\ 
\cQ(5,4)=&\{(1,1,0)\},\\ 
\cQ(5,3)=&\{(3,2,1,0,0),(3,1,1,1,0),(2,2,2,0,0),(2,2,1,1,0),(2,1,1,2,0),\\
&(2,1,1,1,1),(1,1,1,3,0),(1,1,1,2,1)\}.\\
\end{align*}

\noindent As we shall see later, it is not necessary to know the structure of $\cQ(5,2)$ explicitly.

Our analysis is organised according to the following rules: All \textit{visible} colours in $M$ (those represented by Greek alphabet letters) are 2-colours, and both copies of a visible colour are present in $M$. Asterisk entries in $M$ represent $3+$colours. Some asterisk entries appear in $M$ by definition, \textit{e.g.}, each asterisk entry in the first row of $M$ occupies a position in $\{1\}\times[r_2(1)+1,7]$. Another reason why an asterisk entry appears in $M$ is that, if the corresponding position is occupied by a 2-colour $\lambda$, then putting another copy of $\lambda$ to a \textit{free} position (\textit{i.e.}, one that is not occupied by a visible colour) in any proper way (so that the resulting partial vertex colouring $f'$ of $K_6\square K_7$ is proper) leads to a situation, in which no continuation of $f'$ to a proper complete vertex colouring of $K_6\square K_7$ is possible, because at least one pair $\{\lambda,\mu\}$, where $\mu$ is a visible colour, is not good.

If $Q=(3,2,1)$, then (w) $(M)_{4,5}=\zeta$, $(M)_{5,5}=\eta$ and $(M)_{6,6}=\varepsilon$, hence the set $\{\varepsilon,\zeta\}$ is of the type $(2^11^2,2^2)$, which contradicts Claim~\ref{c22}. 

In the case $Q=(2,2,2)$ we are (w) in the situation of Fig.~4. If a 2-colour $\mu$ occupies a position in $\{k\}\times[2l-1,2l]$ for some $k\in[4,6]$ and $l\in[1,2]$, then $\mu=(M)_{4-l,h(k)}$, where $h(k)=\frac12(3k^2-31k+90)$, and $\nu\in C_{3+}$ for each colour $\nu$ occupying a position in $([4,6]\setminus\{k\})\times[5-2l,6-2l]$ (with $\nu\in C_2$ the pair $\{\mu,\nu\}$ is not good). As a consequence, it is easy to see that at least nine positions in $[4,6]\times[1,4]$ are occupied by 3+colours. Besides that, if $\mu=(M)_{i,j}\in C_2$ with $(i,j)\in\{(2,3),(2,4),(3,1),(3,2)\}$, the second copy of $\mu$ must occupy one of the positions $(4,7),(5,5),(6,6)$. Altogether we have $\Sigma\ge3+9+1=13$.

\[
\begin{pmatrix}
\alpha &\beta &\gamma &\delta &\varepsilon &\zeta &\eta\\
\beta &\alpha &. &. &. &. &.\\
. &. &\delta &\gamma &. &. &.\\
. &. &. &. &\zeta &* &.\\
. &. &. &. &. &\eta &*\\
. &. &. &. &* &. &\varepsilon
\end{pmatrix}
\begin{pmatrix}
\alpha &\beta &\gamma &\delta &\varepsilon &\zeta &*\\
\beta &\alpha &. &. &. &. &.\\
. &. &\delta &\gamma &. &. &.\\
. &. &. &. &\zeta &\varepsilon &.\\
. &. &. &. &. &. &*\\
. &. &. &. &. &. &*
\end{pmatrix}
\begin{pmatrix}
\alpha &\beta &\gamma &\delta &\varepsilon &\zeta &*\\
\beta &\alpha &. &. &. &. &.\\ 
. &. &\delta &\gamma &. &. &.\\
. &. &. &. &\zeta &* &.\\
. &. &. &. &. &. &\varepsilon\\
. &. &. &. &. &* &*
\end{pmatrix}
\]
\[
\hskip1mm\mathrm{Fig.\ 4}\hskip33mm\mathrm{Fig.\ 5}\hskip33mm\mathrm{Fig.\ 6}
\]

$Q=(0)$, Fig.~5: Because of $c_2(7)\ge6-c_{3+}\ge2$ we suppose (w) $\eta=(M)_{2,7}\in C_2$ so that $\eta$ is in $\{(3,5),(3,6),(4,3),(4,4)\}$. 

Under the assumption $\eta\in\ro(3)$ we have (w) $\eta=(3,5)$. Let $C_2'$ be the set of 2-colours occupying a position in $[5,6]\times[1,6]$; clearly, $c_{3+}\le4$ implies $|C_2^1|\ge6$. If $\mu\in C_2'$, then from the fact that each pair $\{\mu,\nu\}$ with $\nu\in C_2''=\{\alpha,\beta,\gamma,\delta,\eta\}$ is good one easily gets that the second copy of $\mu$ occupies a position in $[2,3]\times[1,6]$. As $g(4,7,C_2'')=0$ and $g(i,j,C_2'')\le5$ provided that $(i,j)\in[2,3]\times[1,6]$ is a dot position, we obtain $\omega=(M)_{4,7}\in C_{3+}$ (notice that $\omega\notin C_2'$), hence $(M)_{3,7}\in C_2$. Then $\exc(\eta)\ge4$, since we can uncolour vertices $(i,j)$ with $i\in[2,3]$ and $(M)_{i,j}\in C_{3+}$ (here we use $r_{3+}(i)\ge1$ and $C_{3+}\subseteq\co(7)$), as well as the vertices $(5,5),(6,5)$ (independently from the frequencies of $(M)_{5,5}$ and $(M)_{6,5}$) without affecting the completeness of the colour class $\eta$ in the resulting partial colouring.

In the case $\eta\in\ro(4)$ we obtain a contradiction similarly as above.

The assumption $Q=(2,2,0)$ means that (w) $(M)_{4,5}=\zeta$ and $(M)_{5,6}=\varepsilon$, hence the type of the set $\{\varepsilon,\zeta\}$ is $(2^11^2,2^2)$ in contradiction to Claim~\ref{c22}.

For $Q=(2,1,1)$ the situation is (w) depicted on Fig.~6. If $\lambda=(M)_{6,j}\in C_2$, where $j\in[2k-1,2k]$ with $k\in[1,2]$, then $\lambda=(M)_{4-k,5}$. As a consequence, $r_{3+}(6)=4$, $\eta=(M)_{6,5}\in C_2$, and with $\mu=(M)_{2,5}$, $\nu=(M)_{3,5}$ we have $\{\mu,\nu\}\subseteq\ro(6)$. Then, however, $\exc(\eta)\ge3$ ($\mu,\nu$ and at least one 3+colour contribute to the excess of $\eta$).

$Q=(1,1,2)$, Fig.~7: Similarly as above we see that $r_{3+}(6)=4$, $\eta=(M)_{6,7}\in C_2$, $\{(M)_{2,7},(M)_{3,7}\}\subseteq\ro_2(6)$, and so $\exc(\eta)\ge4$.

If $Q=(3,3,1,1,0)$, then (w) $(M)_{3,3}=\delta,(M)_{4,3}=\varepsilon$ and $\gamma\in\{(M)_{5,4},$ $(M)_{6,4}\}$ so that the type of the set $\{\gamma,\delta\}$ contradicts Claim~\ref{c22}.

\[
\begin{pmatrix}
\alpha &\beta &\gamma &\delta &\varepsilon &\zeta &*\\
\beta &\alpha &. &. &. &. &.\\ 
. &. &\delta &\gamma &. &. &.\\
. &. &. &. &* &. &\varepsilon\\
. &. &. &. &. &* &\zeta\\
. &. &. &. &* &* &.
\end{pmatrix}
\begin{pmatrix}
\alpha &\beta &\gamma &\delta &\varepsilon &\zeta &*\\
\beta &\alpha &. &. &. &. &.\\
. &. &\delta &. &. &. &*\\
. &. &\zeta &. &. &. &*\\
. &. &. &\varepsilon &. &. &.\\
. &. &. &. &\gamma &. &*
\end{pmatrix}
\begin{pmatrix}
\alpha &\beta &\gamma &\delta &\varepsilon &\zeta &*\\
\beta &\alpha &. &. &\bullet &\bullet &\bullet\\ 
\bullet &\bullet &\delta &. &* &\bullet &*\\
\bullet &\bullet &\varepsilon &. &* &\bullet &\bullet\\
\bullet &\bullet &. &\zeta &\bullet &\bullet &\bullet\\
\bullet &\bullet &. &. &\bullet &\bullet &\gamma
\end{pmatrix}
\]
\[
\hskip1mm\mathrm{Fig.\ 7}\hskip33mm\mathrm{Fig.\ 8}\hskip33mm\mathrm{Fig.\ 9}
\]

If $Q=(3,2,2,1,0)$, then, having in mind Claim~\ref{c22}, we are (w) in the situation of Fig.~8.  Further, $\eta=(M)_{2,7}\in C_2$ and $\vartheta=(M)_{5,7}\in C_2$, which implies $\eta=(M)_{5,3}$ and $\vartheta=(M)_{2,3}$. Consequently, it is easy to see that both positions in $\{(3,4),(4,6)\}$ are occupied by 3+colours, and, provided that $\mu=(M)_{i,j}\in C_2$ for some $(i,j)\in[3,6]\times[1,2]$, then $(i,j)\in\{(5,1),(5,2)\}$ and $\mu=(M)_{6,3}$. Therefore, $\Sigma\ge4+2+7+r_{3+}(2)\ge14$.

$Q=(3,2,1,1,1)$, (w) Fig.~9 (using Claim~\ref{c22} again): If a bullet position is occupied by a 2-colour $\mu$, then the second copy of $\mu$ occupies a dot position. Therefore, $\Sigma\ge4+(19-6)=17$.

$Q=(3,1,1,1,2)$, (w) Fig.~10: Analogously as in the case of Fig.~9 we obtain $\Sigma\ge5+(19-6)=18$.

Under the assumption $Q=(2,2,2,2,0)$ we have $g(i,7,\{\alpha,\beta,\gamma,\delta,\varepsilon,\zeta\})=1$ for any $i\in[3,6]$ and $g(k,l,\{\alpha,\beta,\gamma,\delta,\varepsilon,\zeta\})\le4$ for any  position $(k,l)\in[2,6]\times[1,6]$ occupied by a colour of $C\setminus\{\alpha,\beta,\gamma,\delta,\varepsilon,\zeta\}$, hence $c_{3+}(7)\ge5$, a contradiction.

\[
\begin{pmatrix}
\alpha &\beta &\gamma &\delta &\varepsilon &\zeta &*\\
\beta &\alpha &. &\bullet &\bullet &\bullet &.\\ 
\bullet &\bullet &\delta &* &\bullet &\bullet &.\\
\bullet &\bullet &\varepsilon &\bullet &* &\bullet &.\\
\bullet &\bullet &. &* &* &\bullet &\gamma\\
\bullet &\bullet &. &\bullet &\bullet &\bullet &\zeta
\end{pmatrix}
\begin{pmatrix}
\alpha &\beta &\gamma &\delta &\varepsilon &\zeta &*\\
\beta &\alpha &. &. &. &. &.\\ 
. &. &\delta &\bullet &. &\bullet &\bullet\\
. &. &. &\varepsilon &\bullet &\bullet &\bullet\\
. &. &\bullet &. &\gamma &\bullet &\bullet\\
. &. &. &. &. &* &\zeta
\end{pmatrix}
\begin{pmatrix}
\alpha &\beta &\gamma &\delta &\varepsilon &\zeta &*\\
\beta &\alpha &. &. &. &. &.\\  
. &. &\delta &* &. &* &.\\
. &. &. &\varepsilon &. &. &*\\
. &. &. &. &\zeta &* &.\\
. &. &. &. &. &. &\gamma
\end{pmatrix}
\]
\[
\hskip1mm\mathrm{Fig.\ 10}\hskip33mm\mathrm{Fig.\ 11}\hskip33mm\mathrm{Fig.\ 12}
\]

If $Q=(2,2,2,1,1)$, because of Claim~\ref{c22} (w) there are two possibilities for the structure of $M$, see Figs.~11 and 12. 

In the case of Fig.~11 a bullet position can be occupied by a 2-colour only if the second copy of that colour appears in $\{2\}\times[3,5]$. A position in $\{(2,6),(2,7),(6,1),(6,2)\}$ is occupied by a 2-colour only if the second copy of that colour is in $\{(3,5),(4,3),(5,4)\}$. Further, it is easy to see that at most one of two colours in $\{i\}\times[1,2]$ with $i\in[3,5]$ is a 2-colour (which is in $\{6\}\times[3,5]$). Therefore $\Sigma\ge2+(9-3)+(4-3)+3\cdot1+r_{3+}(2)\ge13$.

In the situation of Fig.~12 let $k=\max(i\in\{2,3,5\}:(M)_{i,7}\in C_2)$ and $\eta=(M)_{k,7}$. The assumption $k=2$ implies $\eta\in\{(M)_{3,5},(M)_{5,4}\}$.
 
\[
\begin{pmatrix}
\alpha &\beta &\gamma &\delta &\varepsilon &\zeta &*\\
\beta &\alpha &. &. &. &. &\eta\\ 
. &. &\delta &* &\eta &* &*\\
* &* &. &\varepsilon &. &. &*\\
* &* &. &. &\zeta &* &*\\
* &* &. &. &. &. &\gamma
\end{pmatrix}
\begin{pmatrix}
\alpha &\beta &\gamma &\delta &\varepsilon &\zeta &*\\
\beta &\alpha &. &. &. &. &\eta\\ 
* &* &\delta &* &. &* &*\\
. &. &. &\varepsilon &. &. &*\\
. &. &. &\eta &\zeta &* &*\\
* &* &. &. &. &. &\gamma
\end{pmatrix}
\begin{pmatrix}
\alpha &\beta &\gamma &\delta &\varepsilon &\zeta &*\\
\beta &\alpha &. &. &\eta &. &.\\  
. &. &\delta &* &. &* &\eta\\
* &* &. &\varepsilon &. &. &*\\
* &* &. &. &\zeta &* &*\\
. &. &* &. &. &. &\gamma
\end{pmatrix}
\]
\[
\hskip1mm\mathrm{Fig.\ 13}\hskip33mm\mathrm{Fig.\ 14}\hskip33mm\mathrm{Fig.\ 15}
\]

If $\eta=(M)_{3,5}$, see Fig.~13, then $\Sigma\ge13+r_{3+}(2)\ge14$. 

In the case $\eta=(M)_{5,4}$ depicted in Fig.~14 we have $r_{3+}(3)\ge5$. 

If $k=3$ (Fig.~15), then $\eta=(M)_{2,5}$, and from $r_{3+}(2)\ge1$ it follows that $\exc(\alpha)\ge3$, a contradiction.

Fig.~16 corresponds to $k=5$, requiring $\eta=(M)_{2,4}$. If $\vartheta\in C_2$ is in $\{4\}\times[1,2]$, then $\vartheta=(M)_{5,3}$, and, if $\iota\in C_2$ is in $\{6\}\times[1,2]$, then $\iota=(M)_{5,4}$. So, $c_{3+}(1)+c_{3+}(2)\ge4$, which implies $\exc(\alpha)\ge3$.

\[
\begin{pmatrix}
\alpha &\beta &\gamma &\delta &\varepsilon &\zeta &*\\
\beta &\alpha &. &\eta &. &. &.\\  
* &* &\delta &* &. &* &.\\
. &. &. &\varepsilon &. &. &*\\
. &. &. &. &\zeta &* &\eta\\
. &. &. &. &. &* &\gamma
\end{pmatrix}
\begin{pmatrix}
\alpha &\beta &\gamma &\delta &\varepsilon &\zeta &*\\
\beta &\alpha &. &. &. &. &.\\ 
. &. &\delta &. &. &* &.\\
. &. &. &\varepsilon &* &. &.\\
. &. &* &. &* &. &\gamma\\
. &. &. &. &. &* &\zeta
\end{pmatrix}
\begin{pmatrix}
\alpha &\beta &\gamma &\delta &\varepsilon &\zeta &*\\
\beta &\alpha &. &\bullet &\bullet &\bullet &.\\  
\bullet &\bullet &\delta &\bullet &\bullet &\bullet &.\\
\bullet &\bullet &. &\bullet &* &* &\gamma\\
\bullet &\bullet &. &\bullet &* &\bullet &\varepsilon\\
\bullet &\bullet &. &\bullet &\bullet &* &\zeta
\end{pmatrix}
\]
\[
\hskip1mm\mathrm{Fig.\ 16}\hskip33mm\mathrm{Fig.\ 17}\hskip33mm\mathrm{Fig.\ 18}
\]

In the case $Q=(2,2,1,1,2)$, using Claim~\ref{c22}, (w) the description by Fig.~17 applies. Claim~\ref{3C2} implies that a 2-colour occupying a position in $[3,6]\times[1,2]$ does not appear in $\{2\}\times[3,7]$. Therefore, it is easy to see that for any $i\in[3,6]$ at most one of the positions in $\{i\}\times[1,2]$ is occupied by a 2-colour; as a consequence of $c_{3+}\le4$ and $r_{3+}(2)\ge1$ then $\exc(\alpha)\ge3$.

If $Q=(2,1,1,1,3)$, then (w) we have the situation of Fig.~18 with  $\Sigma\ge5+(19-6)=18$ (reasoning as in Fig.~9).

\[
\begin{pmatrix}
\alpha &\beta &\gamma &\delta &\varepsilon &* &*\\
\beta &\alpha &. &. &. &. &.\\ 
. &. &\delta &\gamma &. &. &.\\
. &. &. &. &. &\varepsilon &.\\
. &. &. &. &. &. &*\\
. &. &. &. &. &. &*
\end{pmatrix}
\begin{pmatrix}
\alpha &\beta &\gamma &\delta &\varepsilon &* &*\\
\beta &\alpha &. &. &. &. &.\\ 
. &. &\delta &\bullet &. &. &\bullet\\
. &. &\varepsilon &. &\bullet &. &\bullet\\
. &. &. &. &. &\gamma &\bullet\\
. &. &. &\bullet &\bullet &\bullet &\bullet
\end{pmatrix}
\begin{pmatrix}
\alpha &\beta &\gamma &\delta &\varepsilon &* &*\\
\beta &\alpha &. &. &. &\bullet &\bullet\\ 
. &. &\delta &. &. &. &.\\
. &. &. &\varepsilon &. &. &.\\
. &. &. &. &\gamma &. &.\\
\bullet &\bullet &. &. &. &* &*
\end{pmatrix}
\]
\[
\hskip1mm\mathrm{Fig.\ 19}\hskip33mm\mathrm{Fig.\ 20}\hskip33mm\mathrm{Fig.\ 21}
\]

Clearly, if $r_2(1)=5$, then $r_2(i)=5$ and $r_{3+}(i)=2$ for each $i\in[1,6]$, hence $c_2=\frac12\sum_{i=1}^6r_2(i)=15$, and, by Claims~\ref{eleven}.5, \ref{eleven}.6, $c_{3+}=c_3=4$.

If $Q=(1,1,0)$, we are (w) in the situation of Fig.~19. Each colour of $\co_2(7)$ has its second copy in $[2,4]\times[1,6]$, hence at least $5+2c_2(7)+\sum_{i=2}^4r_{3+}(i)=11+2c_2(7)\ge15$ positions in $[2,4]\times[1,7]$ are occupied by colours of $\{\alpha,\beta,\gamma,\delta,\varepsilon\}\cup\co_2(7)\cup C_{3+}$. Since a colour in $\ro_2(5)\cup\ro_2(6)$ has its second copy in $[2,4]\times[1,6]$, we have $r_2(5)+r_2(6)\le18-(15-3)=6$ and $4=r_{3+}(5)+r_{3+}(6)=14-[r_2(5)+r_2(6)]\ge8$, a contradiction. 

If $Q=(3,2,1,0,0)$, then the set $\{\gamma,\delta\}\subseteq C_2$ is of the type $(2^11^2,2^2)$, which contradicts Claim~\ref{c22}.

$Q=(3,1,1,1,0)$, (w) Fig.~20: A bullet position can be occupied by a colour $\mu\in C_2$ only if $\mu=(M)_{2,3}$. That is why $r_{3+}(6)\ge3$, a contradiction.

$Q=(2,2,2,0,0)$, (w) Fig.~21: If a bullet position is occupied by a colour $\mu\in C_2$, then $\mu\in\{(M)_{3,5},(M)_{4,3},(M)_{5,4}\}$. One can easily see that if $i\in[3,5]$, then at most one of colours in $\{i\}\times[6,7]$ is a 2-colour. Therefore, if $(M)_{2,j}\in C_{3+}$ for both $j=6,7$, then $c_{3+}(6)+c_{3+}(7)\ge3\cdot2+3\cdot1=9$, and there is $j\in[6,7]$ with $c_{3+}(j)\ge5$, a contradiction. Thus, there is $j\in[6,7]$ with $(M)_{2,j}\in C_2$. Then, however, since $(M)_{6,1},(M)_{6,2}\in C_2$ (a consequence of $r_{3+}(6)=2$), the pair $\{(M)_{2,j},(M)_{6,l}\}$ is not good for $l=1,2$.

If $Q=(2,2,1,1,0)$, then (w), by Claim~\ref{c22}, the situation is depicted in Fig.~22. If a 2-colour $\mu$ is in $\{(2,7),(6,1),(6,2)\}$, then $\mu\in\{(M)_{4,3},(M)_{5,4}\}$, and if a 2-colour $\nu$ is in $\{(5,7),(6,6)\}$, then $\nu=(M)_{2,4}$. From $r_{3+}(6)=2$ it follows that there is a 2-colour $\zeta$ in $\{6\}\times[1,2]$; as a consequence then $\omega=(M)_{2,7}\in C_{3+}$ (with $\omega\in C_2$ the pair $\{\omega,\zeta\}$ is not good), $\eta=(M)_{5,7}=(M)_{2,4}\in C_2$, $(M)_{6,6}\in C_{3+}$, and each colour, occupying a position in $\{6\}\times[1,2]$, is a 2-colour. In such a case, however, with $\vartheta=(M)_{4,3}\in\{(M)_{6,1},(M)_{6,2}\}$ the pair $\{\vartheta,\eta\}$ is not good.

$Q=(2,1,1,2,0)$, (w) Fig.~23: If $\zeta\in\{(M)_{3,7},(M)_{6,4}\}\cap C_2$, then $\zeta=(M)_{2,6}$, and if $\eta\in\{(M)_{5,7},(M)_{6,5}\}\cap C_2$, then $\eta=(M)_{2,3}$. Therefore, at least two positions in $\{(3,7),(5,7),(6,4),(6,5)\}$ are occupied by 3+colours. Since $c_{3+}(7)\le4$, at most one position in $\{(3,7),(5,7)\}$ and at least one position in $\{(6,4),(6,5)\}$ is occupied by a 3+colour. Further, from $r_{3+}(6)=2$ it follows that exactly one position in $\{(6,4),(6,5)\}$ and in $\{(3,7),(5,7)\}$ as well is occupied by a 3+colour. Consequently, by Claim~\ref{r23}, $(M)_{2,7}$, $(M)_{6,1}$ and $(M)_{6,2}$ are three distinct 2-colours; this, however, leads to a contradiction, because if $\vartheta\in\{(M)_{2,7},(M)_{6,1},(M)_{6,2}\}\cap C_2$, then necessarily $\vartheta\in\{(M)_{3,6},(M)_{5,3}\}$.

If $Q=(1,1,1,3,0)$, we have (w) $\{\gamma,\delta,\varepsilon\}\cap\ro(6)=\emptyset$. If a position in $\{6\}\times([1,5]\cup\{7\})$ is occupied by a 2-colour $\zeta$, then $\zeta=(M)_{2,6}$, which yields $r_{3+}(6)\ge5$, a contradiction.

\[
\begin{pmatrix}
\alpha &\beta &\gamma &\delta &\varepsilon &* &*\\
\beta &\alpha &. &. &. &. &.\\ 
. &. &\delta &. &. &. &*\\
. &. &. &\varepsilon &. &. &.\\
. &. &. &. &. &\gamma &.\\
. &. &. &. &. &. &*
\end{pmatrix}
\begin{pmatrix}
\alpha &\beta &\gamma &\delta &\varepsilon &* &*\\
\beta &\alpha &. &. &. &. &.\\ 
. &. &\delta &. &. &. &.\\
. &. &. &. &. &\gamma &*\\
. &. &. &. &. &\varepsilon &.\\
. &. &. &. &. &. &*
\end{pmatrix}
\begin{pmatrix}
\alpha &\beta &\gamma &\delta &\varepsilon &* &*\\
\beta &\alpha &. &. &. &. &.\\ 
. &. &* &. &. &\gamma &.\\
. &. &. &* &. &\delta &.\\
. &. &. &. &. &. &\varepsilon\\
. &. &* &* &. &. &.
\end{pmatrix}
\]
\[
\hskip1mm\mathrm{Fig.\ 22}\hskip33mm\mathrm{Fig.\ 23}\hskip33mm\mathrm{Fig.\ 24}
\]

If $Q=(1,1,1,2,1)$, the situation is (w) described by Fig.~24. If a 2-colour $\zeta$ is in $\{6\}\times[1,2]$, then $\zeta=(M)_{5,6}$, hence $r_{3+}(6)\ge3$.

If $Q\in\cQ(5,2)$, we have $\sum_{j=1}^7 q_j=10$. Let $J=\{j\in[1,7]:q_j\ge2\}$. If $|J|\le3$, realise that any colour $\zeta\in C_2\setminus\{\alpha,\beta,\gamma,\delta,\varepsilon\}$ requires existence of a \textit{sufficient} pair $(i,j)\in[2,6]\times[1,7]$, \textit{i.e.}, one satisfying $g(i,j,\{\alpha,\beta,\gamma,\delta,\varepsilon\})\ge3$. If $(i,j)$ is a sufficient pair, then necessarily $j\in J$. Moreover, given $j\in J$, the number of sufficient pairs $(i,j)$ is at most three. This is certainly true if $q_j=3$. On the other hand, if $q_j=2$ and $(M)_{k,l}=(M)_{1,j}$ with $k\ne1$, then, by Claim~\ref{r22} and the fact that $p=2$, $(M)_{k,j}\notin\{\alpha,\beta,\gamma,\delta,\varepsilon\}$, which means that $g(k,j,\{\alpha,\beta,\gamma,\delta,\varepsilon\})=2$, and there are at most three $i$'s such that the pair $(i,j)$ is sufficient. Therefore, $c_2\le5+3\cdot3=14$, which contradicts Claim~\ref{eleven}.3.

So, we have $|J|\ge4$. If $q_1=3$, then $10=\sum_{j=1}^7 q_j\ge3+3\cdot2+1\cdot1=10$, hence $q_2=q_3=q_4=2$, $q_5=1$ and $q_6=q_7=0$. If $\zeta=(M)_{i,j}\in C_2$ with $(i,j)\in[2,6]\times[6,7]$, then, since $g(i,j,\{\alpha,\beta,\gamma,\delta,\varepsilon\})=1$, we have $\zeta\in\co_2(1)\setminus\{\alpha,\beta,\gamma,\delta,\varepsilon\}$. Thus $c_{3+}(6)+c_{3+}(7)\ge2+(10-3)=9$, and there is $j\in[6,7]$ with $c_{3+}(j)\ge5$, a contradiction.

If $q_1\le2$, then $g(i,j,\{\alpha,\beta,\gamma,\delta,\varepsilon\})\le q_j+1\le q_1+1\le3$ for every $(i,j)\in[2,6]\times[1,5]$, hence $g(k,l,\{\alpha,\beta,\gamma,\delta,\varepsilon\})\ge2$ whenever $(k,l)\in[2,6]\times[6,7]$ and $(M)_{k,l}\in C_2$, which implies $q_l\ge1$, $l=6,7$. As a consequence then $10=\sum_{j=1}^7 q_j\ge2|J|+(7-|J|)=|J|+7\ge11$, a final contradiction proving Theorem~\ref{main}. $\qed$
\vskip3mm

\noindent{\large{\bf Acknowledgements.}} This work was supported by the Slovak Research and Development Agency under the contract APVV-19-0153.

\end{document}